# GENERAL DEFINITIONS OF CHAOS
# for CONTINUOUS and DISCRETE-TIME PROCESSES

**By Andrei Vieru**

## *Abstract*

A precise definition of chaos for discrete processes based on iteration already exists. We'll first reformulate it in a more general frame, taking into account the fact that discrete chaotic behavior is neither necessarily based on iteration nor strictly related to compact metric spaces and to bounded functions. Then we'll apply the central idea of this definition to continuous processes. We'll try to see what chaos *is*, regardless of the way it is generated.

## *Introduction*

This paper is motivated by recurrent complaint about the lack of a generally accepted general definition of chaos. One can read on the well known website mathworld.wolfram.com (http://mathworld.wolfram.com/Chaos.html) that many mathematicians, when asked 'what is chaos?', often quote their colleagues, avoiding thus to express their own point of view.

Let's start with Bau-Sen Du's article *'On the nature of chaos'* (arXiv:math.DS/0602585 v1 26 Feb 2006). We find in this article the following definition of chaos:

'We believe that chaos should involve not only nearby points can diverge apart but also faraway points can get close to each other. Therefore, we propose to call a continuous map $f$ from an infinite compact metric space $(X, d)$ to itself chaotic if there exists a positive number $\lambda$ such that for any point $x$ and any nonempty open set $V$ (not necessarily an open neighborhood of $x$) in X there is a point $y$ in V such that $\lim \sup_{n \to \infty} d(f^n(x), f^n(y)) \geq \lambda$ and $\lim \inf_{n \to \infty} d(f^n(x), f^n(y)) = 0$.'

Despite our enthusiasm for Bau-Sen Du's article and for its main underlying idea, we'll keep using a somewhat heavier notation. Here are the main reasons to do so:

1) We'll apply some analogous definitions to both continuous and discrete processes, not only to discrete *iterative* processes.

2) In order to try to reach *general* chaos definitions, we'll not specifically speak about ODE, PDE or about iteration, because we want to see what chaos really *is*, regardless of the way it is obtained.

3) Unlike Bau-Sen Du, we'll consider also mappings of a metric space into *another* metric space, not only onto itself[1].

4) Unlike Bau-Sen Du, we don't suppose that functions are bounded or that metric spaces are compact[2].

---

[1] For maps from a metric space to another metric space that generate chaos see *'Generalized iteration, catastrophes, generalized Sharkovsky's ordering'* arXiv:0801.3755 **[math.DS]**

[2] For unbounded functions generating chaos on non compact metric spaces see *'Periodic helixes, unbounded functions, L-iteration'* arXiv:0802.1401 **[math.DS]**

As far as discrete processes are concerned, we'll consider them in a more general frame, including sequences based not only on the *iteration operator*, but on arbitrary operators (e.g. generalized L-, F- and V-iteration operators[3]).

Doing so, we must be completely aware of the shortcut Bau-Sen Du's uses in his formulations: when he writes 'a function *f* is chaotic if …', we should understand 'the family of sequences generated by iterating *f* is chaotic if…'. A function in itself cannot be neither chaotic nor non-chaotic.

Only a family of functions depending on parameters – in particular a family of sequences depending on parameters (depending either on first term or on some other parameters) – may be chaotic or not.

Besides, even if we accept Bau-Sen Du's shortcut (and we accept it as long as we still keep in mind that it *is* a shortcut), there is another shortcut to be pointed to: we should say instead of 'the function *f* is chaotic if…' 'the function *f* is chaotic *with respect to the iteration operator* if...'.

## *Prelude*

Let $(\psi_a)_{a \in \mathbf{\Omega}}$ be a family of continuous functions, mapping a metric space $(\mathbf{\Theta}, d_1)$ into a metric space $(\mathbf{\Lambda}, d_2)$ and depending on a parameter $a$, whose values are called *initial conditions*. Let $(\mathbf{\Omega}, d_0)$ be the metric space of initial conditions. We'll more often assume $\mathbf{\Theta}$ and $\mathbf{\Lambda}$ are either connected sets or everywhere dense parts of connected sets. In discrete processes, we'll only assume $\mathbf{\Theta}$ has at least one cluster point (usually $\infty$). When we'll consider continuous processes, we'll more often assume $(\mathbf{\Theta}, d_1)$ is a *path-connected metric space*. However, our definitions still hold for metric spaces that are not path-connected. The restrictions contained in our definitions drop by themselves if the set of paths in the metric space (or in its considered subsets) is empty.

We'll call the elements of $\mathbf{\Omega}$ *points*, while the elements of the domain $\mathbf{\Theta}$ of the functions $\psi_a$ will be called *moments*, even if $\mathbf{\Theta}$ may not necessarily be a subset of $\mathbf{R}$. In this article, we'll not consider cases in which $\mathbf{\Theta}$ and $\mathbf{\Omega}$ are sets of fractal dimension. Since there is no possible confusion, we'll designate $d_0(x, y)$, $d_1(z, t)$ and $d_2(u, v)$ by $|x - y|$, $|z - t|$ and $|u - v|$.

## 1. Sensitive dependence on initial conditions

1.1. For discrete systems we'll just transcribe the classical formulation of this property: $\exists \lambda \in \mathbf{R}\ \forall x \in \mathbf{\Omega}\ \forall \varepsilon \in \mathbf{R}\ \exists y \in \mathbf{\Omega}\ \exists n$ in $\mathbf{N}$ (i.e. in[4] $\mathbf{\Theta}$) $|x - y| < \varepsilon \wedge |f^n(x) - f^n(y)| > \lambda$

1.2. Or, for continuous systems:
$\exists \lambda \in \mathbf{R}\ \forall x \in \mathbf{\Omega}\ \forall \varepsilon \in \mathbf{R}\ \exists y \in \mathbf{\Omega}\ \exists z$ in $\mathbf{\Theta}$ $|x - y| < \varepsilon \wedge |f_x(z) - f_y(z)| > \lambda$

---

[3] For a definition of generalized L-iteration, see our paper *'Periodic helixes, unbounded Functions, L-iteration' 'iteration'* arXiv:0802.1401 **[math.DS]**. For a definition of generalized F- and V-iteration see our paper *'Generalized Iteration, Catastrophes and Generalized Sharkovsky's ordering'* arXiv:0801.3755 **[math.DS]**

[4] We adopt here the traditional notation. However, instead of $f_n(x)$ – or of $f^n(x)$ if the sequences are generated by iteration – we'll prefer further to designate by $f_x(n)$ the *n*-th term of a sequence from a family of sequences depending on a parameter – which, for example may be its first term – when the value of the parameter is $x$. The notation will thus be analogous to 1.2.

(An attempt to establish the beginning of a classification of the different types of sensitive dependence may be found at the end of this paper, in the Appendix.)

Since sensitive dependence on initial conditions does not imply chaos we need a stronger property, we'll define such a property and we'll call it *chaotic dependence on arbitrarily close or faraway initial conditions*. We'll distinguish two different varieties of this property, namely the disjoint and the cross-graph chaotic dependence:

## 2. 'Disjoint' chaotic dependence on arbitrarily close or faraway initial conditions[5]

2.1.1. Let S be a connected compact subset of $\Theta$ containing at least an open set and let B(S) designate its boundary. Let $\Pi(S_1, S_2)$ designate the set of paths connecting any point $w_1$ in $S_1$ to any point $w_2$ in $S_2$. Let $(\psi_a)_{a\in\Omega}$ be a family of continuous functions that map a metric space $(\Theta, d_1)$ into a metric space $(\Lambda, d_2)$. The family $(\psi_a)_{a\in\Omega}$ displays *strong disjoint chaotic dependence on arbitrarily close or faraway initial conditions* if $\exists\mu>0$ $\forall\alpha\in\Omega$ $\forall\beta\in\Omega$ $|\beta-\alpha|>0 \Rightarrow$

1) $\forall\varepsilon_1\in R$ $\exists S_1\subset\Theta$ $[\forall z\in S_1-B(S_1)$ $0<|\psi_\alpha(z)-\psi_\beta(z)|<\varepsilon_1$ $\land$ $\forall z\in B(S_1)$ $|\psi_\alpha(z)-\psi_\beta(z)|=\varepsilon_1]$

2) For every $\varepsilon$ such as $\exists S\subset\Theta$ $[\forall z\in S-B(S)$ $0<|\psi_\alpha(z)-\psi_\beta(z)|<\varepsilon$ $\land$ $\forall z\in B(S)$ $|\psi_\alpha(z)-\psi_\beta(z)|=\varepsilon]$ there is an $\varepsilon_2<\varepsilon$ $\exists S_2\subset\Theta$ $[S_2\cap S=\varnothing]$ $\land$ $[\forall z\in S_2-B(S_2)$ $0<|\psi_\alpha(z)-\psi_\beta(z)|<\varepsilon_2]$ $\land$ $[\forall z\in B(S_2)$ $|\psi_\alpha(z)-\psi_\beta(z)|=\varepsilon_2]$

3) $\forall P\in\Pi(S, S_2)$ $\exists\gamma\in P$ $|\psi_\alpha(\gamma) - \psi_\beta(\gamma)|>\mu$.

2.1.2. In particular, if $\Theta\subset\mathbf{R}$, the definition may be formulated as follows: The family $(\psi_a)_{a\in\Omega}$ displays *strong disjoint chaotic dependence on arbitrarily close or faraway initial conditions* if $\exists\mu>0$ $\forall\alpha\in\Omega$ $\forall\beta\in\Omega$ $|\beta-\alpha|>0 \Rightarrow$

1) $\forall\varepsilon_1$ $\exists x_1\in\Theta$ $\exists y_1\in\Theta$ $(x_1<y_1)$ $[|\psi_\alpha(x_1)-\psi_\beta(x_1)|=\varepsilon_1]$ $\land$ $[|\psi_\alpha(y_1)-\psi_\beta(y_1)|=\varepsilon_1]$ $\land$ $[\forall z\in(x_1, y_1)$ $0<|\psi_\alpha(z)-\psi_\beta(z)|<\varepsilon_1]$

2) For every $\varepsilon$ such as $\exists x\in\Theta$ $\exists y\in\Theta$ $[x<y]$ $\land$ $[|\psi_\alpha(x)-\psi_\beta(x)|=\varepsilon]$ $\land$ $[|\psi_\alpha(y)-\psi_\beta(y)|=\varepsilon]$ $\land$ $[\forall z\in(x, y)$ $0<|\psi_\alpha(z)-\psi_\beta(z)|<\varepsilon]$ there is an $\varepsilon_2<\varepsilon$ $\exists x_2\in\Theta$ $\exists y_2\in\Theta$ $(x_2<y_2)$ $[[x_2, y_2]\cap[x, y]=\varnothing]$ $\land$ $[|\psi_\alpha(x_2)-\psi_\beta(x_2)|=\varepsilon_2]$ $\land$ $[\forall z\in(x_2, y_2)$ $0<|\psi_\alpha(z)-\psi_\beta(z)|<\varepsilon_2]$

3) $\exists w\in[\min(y, y_2), \max(x, x_2)]$ $|\psi_\alpha(w)-\psi_\beta(w)|>\mu$

2.2.1. Again, let S be a connected compact subset of $\Theta$ containing at least an open set and let B(S) designate its boundary. Let $\Pi(S_1, S_2)$ designate the set of paths connecting any point $w_1$ in $S_1$ to any point $w_2$ in $S_2$. Let $V_\chi$ be the set of neighborhoods of $\chi$. Let $(\psi_a)_{a\in\Omega}$ a family of continuous functions that map a metric space $(\Theta, d_1)$ into a metric space $(\Lambda, d_2)$. The family $(\psi_a)_{a\in\Omega}$ displays *weak disjoint*

---

[5] this concept is, in its essence, not very complicated: whatever the difference between two initial conditions we'll always find intervals (or, in the general case, sub-domains containing open sets) on which the values of the two functions are as near one to another as we want, without 'touching' or 'crossing' one another. We spent a lot of symbols – so the definition seems heavy – in order to eliminate from its area completely irrelevant cases like that of two straight lines that cross each other in some point. (In the neighborhood of such a point, we find tiny intervals on which the condition would have been satisfied, wouldn't we introduce some restrictions in the definition.) We also wanted to eliminate cases when two functions like, e.g., $ax\sin(1/x)$ meet in some point ('analogous' to 0 in the example).

*chaotic dependence on arbitrarily close or faraway initial conditions* if $\exists\mu>0$ $\forall\alpha\in\Omega$ $\forall\chi\in\Omega$ $\forall v\in V_\chi$ $\exists\beta\in v$ ($\beta\neq\alpha$) such as

1) $\forall\varepsilon_1\in R$ $\exists S_1\subset\Theta$ $[\forall z\in S_1-B(S_1)$ $0<|\psi_\alpha(z)-\psi_\beta(z)|<\varepsilon_1$ $\wedge$ $\forall z\in B(S_1)$ $|\psi_\alpha(z)-\psi_\beta(z)|=\varepsilon_1]$

2) For every $\varepsilon$ such as $\exists S\subset\Theta$ $[\forall z\in S-B(S)$ $0<|\psi_\alpha(z)-\psi_\beta(z)|<\varepsilon$ $\wedge$ $\forall z\in B(S)$ $|\psi_\alpha(z)-\psi_\beta(z)|=\varepsilon]$ there is an $\varepsilon_2<\varepsilon$ $\exists S_2\subset\Theta$ $[S_2\cap S=\varnothing]$ $\wedge$ $[\forall z\in S_2-B(S_2)$ $0<|\psi_\alpha(z)-\psi_\beta(z)|<\varepsilon_2]$ $\wedge$ $[\forall z\in B(S_2)$ $|\psi_\alpha(z)-\psi_\beta(z)|=\varepsilon_2]$

3) $\forall P\in\Pi(S, S_2)$ $\exists\gamma\in P$ $|\psi_\alpha(\gamma)-\psi_\beta(\gamma)|>\mu$.

2.2.2. In particular let $\Theta\subset R$. Let again $V_\chi$ be the set of neighborhoods of $\chi$. The family $(\psi_a)_{a\in\Omega}$ displays *weak disjoint chaotic dependence on arbitrarily close or faraway initial conditions* if $\exists\mu>0$ $\forall\alpha\in\Omega$ $\forall\chi\in\Omega$ $\forall v\in V_\chi$ $\exists\beta\in v$ ($\beta\neq\alpha$) such as

1) $\forall\varepsilon_1$ $\exists x_1\in\Theta$ $\exists y_1\in\Theta$ $(x_1<y_1)$ $[|\psi_\alpha(x_1)-\psi_\beta(x_1)|=\varepsilon_1]$ $\wedge$ $[|\psi_\alpha(y_1)-\psi_\beta(y_1)|=\varepsilon_1]$ $\wedge$ $[\forall z\in(x_1, y_1)$ $0<|\psi_\alpha(z)-\psi_\beta(z)|<\varepsilon_1]$

2) For every $\varepsilon$ such as $\exists x\in\Theta$ $\exists y\in\Theta$ $[x<y]$ $\wedge$ $[|\psi_\alpha(x)-\psi_\beta(x)|=\varepsilon]$ $\wedge$ $[|\psi_\alpha(y)-\psi_\beta(y)|=\varepsilon]$ $\wedge$ $[\forall z\in(x, y)$ $0<|\psi_\alpha(z)-\psi_\beta(z)|<\varepsilon]$ there is an $\varepsilon_2<\varepsilon$ $\exists x_2\in\Theta$ $\exists y_2\in\Theta$ $(x_2<y_2)$ $[[x_2, y_2]\cap[x, y]=\varnothing]$ $\wedge$ $[|\psi_\alpha(x_2)-\psi_\beta(x_2)|=\varepsilon_2]$ $\wedge$ $[|\psi_\alpha(y_2)-\psi_\beta(y_2)|=\varepsilon_2]$ $\wedge$ $[\forall z\in(x_2, y_2)$ $0<|\psi_\alpha(z)-\psi_\beta(z)|<\varepsilon_2]$

3) $\exists w\in[\min(y, y_2), \max(x, x_2)]$ $|\psi_\alpha(w)-\psi_\beta(w)|>\mu$

# 3. 'Cross-graph' chaotic dependence on arbitrarily close or faraway initial conditions

3.1.1. We'll use the same symbols as in the precedent chapter. The family $(\psi_a)_{a\in\Omega}$ shows *strong cross-graph chaotic dependence on arbitrarily close or faraway initial conditions* if $\exists\mu>0$ $\forall\alpha\in\Omega$ $\forall\beta\in\Omega$ $|\beta-\alpha|>0\Rightarrow$

1) $\forall\varepsilon_1\in R$ $\exists S_1\subset\Theta$ $[\forall z\in S_1-B(S_1)$ $0<|\psi_\alpha(z)-\psi_\beta(z)|<\varepsilon_1$ $\wedge$ $\forall z\in B(S_1)$ $\psi_\alpha(z)=\psi_\beta(z)]$

2) for every $\varepsilon$ such as $\exists S\subset\Theta$ $[\forall z\in S-B(S)$ $0<|\psi_\alpha(z)-\psi_\beta(z)|<\varepsilon$ $\wedge$ $\forall z\in B(S)$ $\psi_\alpha(z)=\psi_\beta(z)]$ there is an $\varepsilon_2<\varepsilon$ $\exists S_2\subset\Theta$ $[S_2\cap S=\varnothing]$ $\wedge$ $[\forall z\in S_2-B(S_2)$ $0<|\psi_\alpha(z)-\psi_\beta(z)|<\varepsilon_2]$ $\wedge$ $[\forall z\in B(S_2)$ $\psi_\alpha(z)=\psi_\beta(z)]$

3) $\forall P\in\Pi(S, S_2)$ $\exists\gamma\in P$ $|\psi_\alpha(\gamma)-\psi_\beta(\gamma)|>\mu$.

3.1.2. In particular, if $\Theta\subset R$, the definition may be formulated as follows: *cross-graph chaotic dependence on arbitrarily close or faraway initial conditions* if $\exists\mu>0$ $\forall\alpha\in\Omega$ $\forall\beta\in\Omega$ $|\beta-\alpha|>0\Rightarrow$

1) $\forall\varepsilon_1$ $\exists x_1\in\Theta$ $\exists y_1\in\Theta$ $(x_1<y_1)$ $[|\psi_\alpha(x_1)=\psi_\beta(x_1)|]$ $\wedge$ $[|\psi_\alpha(y_1)=\psi_\beta(y_1)|]$ $\wedge$ $[\forall z\in(x_1, y_1)$ $0<|\psi_\alpha(z)-\psi_\beta(z)|<\varepsilon_1]$

2) For every $\varepsilon$ such as $\exists x\in\Theta$ $\exists y\in\Theta$ $[x<y]$ $\wedge$ $[|\psi_\alpha(x)-\psi_\beta(x)|=0]$ $\wedge$ $[|\psi_\alpha(y)-\psi_\beta(y)|=0]$ $\wedge$ $[\forall z\in(x, y)$ $0<|\psi_\alpha(z)-\psi_\beta(z)|<\varepsilon]$ there is an $\varepsilon_2<\varepsilon$ $\exists x_2\in\Theta$ $\exists y_2\in\Theta$ $(x_2<y_2)$ $[[x_2, y_2]\cap[x, y]=\varnothing]$ $\wedge$ $[|\psi_\alpha(x_2)-\psi_\beta(x_2)|=0]$ $\wedge$ $[|\psi_\alpha(y_2)-\psi_\beta(y_2)|=0]$ $\wedge$ $[\forall z\in(x_2, y_2)$ $0<|\psi_\alpha(z)-\psi_\beta(z)|<\varepsilon_2]$

3) $\exists w\in[\min(y, y_2), \max(x, x_2)]$ $|\psi_\alpha(w)-\psi_\beta(w)|>\mu$

3.2.1. The family $(\psi_a)_{a\in\Omega}$ shows *weak cross-graph chaotic dependence on arbitrarily close or faraway initial conditions* if $\exists\mu>0$ $\forall\alpha\in\Omega$ $\forall\chi\in\Omega$ $\forall v\in V_\chi$ $\exists\beta\in v$ ($\beta\neq\alpha$) such as

1) $\forall\varepsilon_1\in R$ $\exists S_1\subset\Theta$ $[\forall z\in S_1-B(S_1)$ $0<|\psi_\alpha(z)-\psi_\beta(z)|<\varepsilon_1$ $\wedge$ $\forall z\in B(S_1)$ $\psi_\alpha(z)=\psi_\beta(z)]$

2) For every $\varepsilon$ such as $\exists S\subset\Theta$ $[\forall z\in S-B(S)$ $0<|\psi_\alpha(z)-\psi_\beta(z)|<\varepsilon$ $\wedge$ $\forall z\in B(S)$ $\psi_\alpha(z)=\psi_\beta(z)]$ there is an $\varepsilon_2<\varepsilon$ $\exists S_2\subset\Theta$ $[S_2\cap S=\varnothing]$ $\wedge$ $[\forall z\in S_2-B(S_2)$ $0<|\psi_\alpha(z)-\psi_\beta(z)|<\varepsilon_2]$ $\wedge$ $[\forall z\in B(S_2)$ $\psi_\alpha(z)=\psi_\beta(z)]$

3) $\forall P\in\Pi(S, S_2)$ $\exists\gamma\in P$ $|\psi_\alpha(\gamma)-\psi_\beta(\gamma)|>\mu$.

3.2.2. Assume $\Theta \subset \mathbf{R}$. Again, let $V_\chi$ be the set of neighborhoods of $\chi$. The family $(\psi_a)_{a\in\mathbf{\Omega}}$ shows *weak cross-graph chaotic dependence on arbitrarily close or faraway initial conditions* if $\exists\mu>0$ $\forall\alpha\in\mathbf{\Omega}$ $\forall\chi\in\Theta$ $\forall\nu\in V_\chi$ $\exists\beta\in\nu$ ($\beta\neq\alpha$) such as

1) $\forall\varepsilon_1$ $\exists x_1\in\Theta$ $\exists y_1\in\Theta$ $(x_1<y_1)$ $[|\psi_\alpha(x_1)=\psi_\beta(x_1)|]$ $\wedge$ $[|\psi_\alpha(y_1)=\psi_\beta(y_1)|]$ $\wedge$ $[\forall z\in(x_1,\ y_1)$ $0<|\psi_\alpha(z)-\psi_\beta(z)|<\varepsilon_1]$

2) For every $\varepsilon$ such as $\exists x\in\Theta$ $\exists y\in\Theta$ $[x<y]$ $\wedge$ $[|\psi_\alpha(x)-\psi_\beta(x)|=0]$ $\wedge$ $[|\psi_\alpha(y)-\psi_\beta(y)|=0]$ $\wedge$ $[\forall z\in(x,\ y)$ $0<|\psi_\alpha(z)-\psi_\beta(z)|<\varepsilon]$ there is an $\varepsilon_2<\varepsilon$ $\exists x_2\in\Theta$ $\exists y_2\in\Theta$ $(x_2<y_2)$ $[[x_2,\ y_2]\cap[x,\ y]=\varnothing]$ $\wedge$ $[|\psi_\alpha(x_2)=\psi_\beta(x_2)|=0]$ $\wedge$ $[|\psi_\alpha(y_2)=\psi_\beta(y_2)|=0]$ $\wedge$ $[\forall z\in(x_2,\ y_2)$ $0<|\psi_\alpha(z)-\psi_\beta(z)|<\varepsilon_2]$

3) $\exists w\in[\min(y,\ y_2),\ \max(x,\ x_2)]$ $|\psi_\alpha(w)-\psi_\beta(w)|>\mu$

# 4. Chaotic dependence on arbitrarily close or faraway initial conditions in discrete processes

4.0. We don't use here the classical notation for function iteration because we aim at a general chaos definition, regardless of the way chaos is obtained. In our notation, $u_\alpha(n)$ designates the *n*-th term of the sequence (when the parameter value the family of sequences depend on equals $\alpha$). It is necessary to make our definition fit not only with sequences based on iteration of $\Theta\rightarrow\Theta$ maps, but also for arbitrary families of sequences. For instance, with sequences based on F- or V- or L-iteration of $\Theta^n\rightarrow\Theta^m$ (*m*<*n*) maps.

Although for discrete processes (i.e. for families of sequences), the distinction between the disjoint metric independence and the cross-graph metric independence with respect to faraway initial conditions can easily be established, one can doubt about its relevance. We'll not indulge here in such details.

4.1. The family of sequences $(u_a(n))_{a\in\mathbf{\Omega}}$ shows chaotic dependence on (arbitrarily close or faraway) initial conditions if

1)$\forall\alpha$ $\forall\chi\in\mathbf{\Omega}$ $\forall\nu\in V_\chi$ $\exists\beta\in\nu$ $|\beta-\alpha|>0$ $\Rightarrow$ $\forall\varepsilon>0$ $\forall n\in\mathbf{N}$ $\exists m>n$ $[|u_\alpha(m)-u_\beta(m)|<\varepsilon]$

2) $\forall\alpha$ $\forall\chi\in\mathbf{\Omega}$ $\forall\nu\in V_\chi$ $\exists\beta\in\nu$ $|\beta-\alpha|>0$ $\Rightarrow$ $\forall\varepsilon>0$ $\forall n\in\mathbf{N}$ $\exists m>n$ $[|u_\alpha(m)-u_\beta(m)|\geq\lambda]$

4.2. We would like to formulate also this definition in Bau-Du Sen's style: we propose to call a set of (not necessarily everywhere) continuous maps $\{f_1, f_2, ..., f_n\}$ from an infinite (not necessarily compact) metric space $(X, d_1)$ to an infinite (not necessarily compact) metric space $(Y, d_2)$ *'chaotic with respect to* $\mathbf{A}$' if there is an algorithm[6] $\mathbf{A}$ using $\{f_1, f_2, ..., f_n\}$ to generate a strictly deterministic family of sequences $(u_x(n))$ of points in Y depending on their initial terms, a positive number $\lambda$ such that for any point *x* in X and any nonempty open set V (not necessarily an open neighborhood of *x*) in X there is a point y in V such that $\lim\sup_{n\rightarrow\infty}d_2(u_x(n),\ u_y(n))\geq\lambda$ and $\lim\inf_{n\rightarrow\infty}d_2(u_x(n),\ u_y(n))=0$. The set $\{f_1, f_2, ..., f_n\}$ will be simply called 'chaotic' if there is no possible confusion concerning $\mathbf{A}$.

# 5. Chaotic systems of the first kind

5.0. In this chapter we'll drop terms like *strong* or *weak*, offering the reader the possibility to restore them.

---

[6] such an algorithm may be easily conceived using and/or combining generalized F-iteration and/or V-iteration and/or L-iteration.

5.1. A family of functions is a *chaotic system **of the first kind*** if it has cross-graph chaotic dependence on arbitrarily close or faraway initial conditions but not disjoint chaotic dependence on arbitrarily close or faraway initial conditions.

5.2. Example: $\mathbf{\Omega} = [0, 1]$   $\mathbf{\Theta} = ]0, 1]$   $f_a(x) = \{\sin[a\ln(1/2ax)/x]+1\}/2$.

5.3. $\mathbf{\Omega} = \mathbf{R}$, $\mathbf{\Theta} = \mathbf{R}$   $f_a(x) = \sin(ax)$. This chaotic system[7] is of the first kind.

# 6. Chaotic systems of the second kind

6.0. In this chapter we'll drop terms like *strong* or *weak*, offering the reader the possibility to restore them.

6.1. A family of functions is a *chaotic system **of the second kind*** if it has disjoint *chaotic dependence on arbitrarily close or faraway initial conditions*, **but not** *cross-graph chaotic dependence on arbitrarily close or faraway initial conditions*.

# 7. Chaotic systems of the third kind

7.0. In this chapter we'll drop terms like *strong* or *weak*, offering the reader the possibility to restore them.

7.1. A family is a *chaotic system of the third kind* if it has points are chaotic points of the third kind *disjoint chaotic dependence on arbitrarily close or faraway initial conditions* **and** *cross-graph chaotic dependence on arbitrarily close or faraway initial conditions*.

7.2. Example: Let's write $F_a(x) = a(\sin\pi x + \cos2x + x)$ and let's designate by $F_a^n(x)$ the *n*-th iterate of $F_a(x)$. (One should note that $\pi$ is irrational, while 2 is rational,

---

[7] It may seem strange that families of periodic functions might be considered as chaotic. However, if we chose two close incommensurable values – $\alpha$ and $\beta$ – of the parameter $a$ and sufficiently high $x$ values, we'll see that it is as difficult to 'calculate', to 'evaluate' or to 'predict' the value of $\sin(\alpha x)$–$\sin(\beta x)$ as it would be to calculate, evaluate or predict, for arbitrarily close $\alpha$ and $\beta$ values and sufficiently high $n$, the difference $u_\alpha(n) - u_\beta(n)$ when, for instance, $u_\alpha(n)$ and $u_\beta(n)$ designates $f^n(\alpha)$ and $f^n(\beta)$ with $f(x)=3.57x(1 - x)$. To calculate values for big $n$ or $x$ values properly, sensitive dependence on initial conditions makes us need more and more digits. On the other hand, if all values of the parameter $a$ the family $\sin(ax)$ depends on where commensurable, there would have been sensitive dependence on initial conditions here, but there would have been no first kind chaos (and *a fortiori* no chaos at all). As a matter of fact, one can extract from the $(\sin(ax))_{a\in\mathbf{R}}$ chaotic (of the first kind) family, a discrete chaotic family, namely the family of sequences $(u_a(n))_{a\in\mathbf{R}}$, where for every $n$ and for every $a$ we'll write $u_a(n) = \sin(an) - \sin n$. For irrational $a$ values incommensurable with $\pi$, these sequences are aperiodic and unpredictable (in the sense that the greater the value assigned to $n$, more digits you'll need to work with in order to properly predict the $u_a(n)$ value). For $a$ values commensurate with $\pi$, the sequence $u_a(n) = \sin(an) - \sin n$ will be periodic. One could (wrongly, because there is no iteration here in the strict sense) say it constitutes an 'orbit'. This example lets someone think about the analogy between periodic orbits and rationals on one hand, between strange attractors and irrationals on the other hand. (We intuitively perceive as 'chaotic' the coexistence of various periodicities. Eratosthene's sieve – which 'generates' the 'chaotic' (aperiodic) set of primes by means of various periodicities – is but one example. As Vladimir Arnold showed, multiplicative groups in finite fields of order $p^k$ ($k{\geq}2$) constitute another beautiful example. However, as $p{\to}\infty$, Galois fields only *approximate* chaos, as rationals approximate irrationals)

so the function is quasi-periodic. In fact the periodicity of a whole family of functions – not only of the functions in some family, all considered separately – is not incompatible with sensitive dependence on initial conditions **but is incompatible with cross-graph chaotic dependence[8] on arbitrarily close or faraway initial conditions**, **and a fortiori with disjoint chaotic dependence on arbitrarily close or faraway initial conditions**.) Then, for every integer $n \geq 1$, $G_a{}^n(x) = F_a{}^n(x) - a^n x$ ($x \in \mathbf{R}$, $a \in \mathbf{R}$) is a chaotic system of the third kind[9].

7.4. Another example: $f_a(x) = \sum_p [\sin(ax/p)]/p^2$ (where $p$ goes through all primes)

### Conclusion

There is a serious question about any definition or set of conditions. Are they recursive? We mean, given some concrete object, the question is: is there any finite-step algorithm that might provide the proof that the given object satisfies or not the definition or the set of conditions. Actually, we haven't examined seriously this question. We submit it to the reader's attention.

As far as we are concerned, we are rather pessimistic: we think that a general recursive chaos definition simply does not exist. We think that recursive chaos definitions exist only for particular forms of chaos.

### Acknowledgements
We express our deep gratitude to Robert Vinograd for useful discussions.

AndreiVieru

---

[8] Suffice to consider the example of the family $f_a(x) = \sin(ax)$

[9] $F_a{}^n(x) - a^n x = 0$ is the only one equation of the form $F_a{}^n(x) - bx = 0$ that has infinitely many roots

APPENDIX

### A.0. An attempt to classify sensitive dependences

Although this is not the main topic of this print, the definitions of sensitive dependence on initial conditions (see 1.1. and 1.2.) suggest some steps toward a classification of sensitive dependence in terms of non-uniform convergence:

### A.1. Insensitivity to initial conditions

A.1.1. $\alpha \in \Omega$ is a point of *insensitivity to close initial conditions* if, for any sequence of points $u(n) \neq \alpha$ converging to $\alpha$, the sequence of functions $(\psi_{u(n)})$ converges **uniformly** to $\psi_\alpha$ on all $\Theta$.

A.1.2. The family $(\psi_a)_{a \in \Omega}$ is *insensitive to close initial conditions* if all values of parameter $a$ are points of insensitivity to *close* initial conditions.

### A.2. Smooth sensitive dependence on initial conditions

A.2.1. $\alpha \in \Omega$ is a point of *smooth sensitive dependence on (close) initial conditions* if, for *any* sequence of points $u(n) \neq \alpha$ converging to $\alpha$, the sequence of functions $(\psi_{u(n)})$ converges **non-uniformly** to $\psi_\alpha$, on all $\Theta$, save, possibly, on its boundaries.

A.2.2. The family $(\psi_a)_{a \in \Omega}$ displays *smooth sensitive dependence on close initial conditions* if all values of the parameter $a$ are points of smooth sensitive dependence on close initial conditions.

A.2.3. Three examples of smooth sensitive dependence on close initial conditions:

$\Omega = [0, 1]$     $\Theta = \,]0, 1]$     $g_a(x) = \{\sin[a\ln(1/2ax)/x]+1\}/2$

$\Omega = \Theta = \mathbf{R}$            $h_a(x) = ax$

$\Omega = [0, 1]$     $\Theta = [0, 1[$     $f_a(x) = [\sin(2\pi a/(1-x))+1]/2$.

### A.3. Points of discontinuous sensitive dependence

A.3.1. $\alpha \in \Omega$ is a *point of discontinuity*[10] *within a context of sensitive dependence on close initial conditions* if, for *almost* any sequence of points $u(n) \neq \alpha$ converging to $\alpha$, the sequence of functions $(\psi_{u(n)})$ converges non-uniformly on all $\Theta$ – except, possibly, on its boundaries – to some $\psi_\beta \neq \psi_\alpha$.

(We cannot dispense with the word '*almost*'. This can be shown on the following example: $\Omega = [0, 1]$   $\Theta = [0, 1[$   $f_a(x) = [\sin(2\pi a/(1-x))+1]/2$ if $a$ is irrational and $f_a(x) = ax$ if $a$ is rational. If we chose a sequence $u(n)$ of irrationals converging to 0.5 and a sequence $v(n)$ of rationals converging to 0.5, then the sequence of functions $(f_{u(n)}(x))$ will converge non-uniformly to $y = [\sin(\pi/(1-x))+1]/2$, while the sequence of functions $(f_{v(n)}(x))$ will converge uniformly to $y = x/2$. If $\Theta = \mathbf{R} - \{1\}$, then we'll have the same convergences, with the only difference that $(f_{v(n)}(x))$ will converge **non-uniformly** to $y = x/2$.)

A.3.2. The family $(\psi_a)_{a \in \Omega}$ is *discontinuous within a context of sensitive dependence on close initial conditions* if all values of parameter $a$ are either *points of sensitive dependence on close initial conditions* or *points of discontinuity within a context of sensitive dependence on close initial conditions*.

---

[10] We'll not develop here the idea of point of discontinuity, though it might be an interesting concept, even in the context of *insensitive* dependence. See the catastrophes we describe in our paper *'Generalized Iteration, Catastrophes and Generalized Sharkovsky's ordering'*. arXiv:0801.3755 **[math.DS]**

A.4. **Points of disordered sensitive dependence**

A.4.1. $\alpha \in \Omega$ is a *point of total disorder* if, for any sequence of points $u(n) \neq \alpha$ converging to $\alpha$, the sequence of functions $(\psi_{u(n)})$ does not converge on the interior of $\Theta$.

A.4.2. The family $(\psi_a)_{a \in \Omega}$ is *a totally disordered family* if all values of parameter $a$ are *points of total disorder*.

A.4.3. This type of family might not exist, unless we accept to define functions *via* infinite stochastic processes. In order to build an example of a totally disordered family one can take an everywhere discontinuous function $\Xi(x) : (0, 1) \to (0, 1)$ and write $\psi_\alpha(x) = x\Xi(a)$. For example, introducing a probabilistic measure on $\mathbf{N}$, we can construct the everywhere discontinuous function $\Xi(x) : (0, 1) \to (0, 1)$ by choosing for every $x$ the $n$-th randomly chosen $y$ value in the interval for which $\{\sin[1/(1-y)]+1\}/2 = x$.

A.4.4. We'll say, in this weaker variant of definition 1.2.4.1., that $\alpha$ is a *point of dense disorder* if there is a dense set of sequences of points $u(n) \neq \alpha$ converging to $\alpha$, such as the sequence of functions $(\psi_{u(n)})$ does not converge on the interior of $\Theta$.

A.4.5. We'll say that the family $(\psi_a)_{a \in \Omega}$ is dense-disordered if there is a dense set of points of dense disorder (in the sense of 1.2.4.4.).

A.4.6. Example: Let $\Theta$ be $[0, 1[$, let $\Omega$ be $]0, 1]$ and let's first define a function $\vartheta(x)$ as follows:

$\vartheta(x) = 3/4$ if $x$ is a rational

$\vartheta(x) = 1$ if $x$ is transcendental

$\vartheta(x) = 1/n$ if $x$ is an algebraic irrational and if $n$ is the lowest degree of a polynomial one of whose roots is $x$.

We can now define the map family $\psi_\alpha(x) = \{\sin\{[\sin(1/\vartheta(a))+1]/(2-2x)\}+1\}/2$